\input amstex 
\documentstyle{amsppt} 
\magnification=1200 
\NoRunningHeads
\topmatter 
\title 
Smoothness of the uniformization of two-dimensional linear foliation on  
 torus with nonstandard metric\endtitle 
\author A.A.Glutsyuk\endauthor 
\thanks Research supported by CRDF grant RM1-229, by INTAS grant 93-0570-ext, 
by Russian Foundation for Basic Research (RFBR) grant 98-01-00455, by State 
Scientific Fellowship of Russian Academy of Sciences for young scientists, 
by the European Post-Doctoral Institute joint fellowship of Max-Planck Institut f\"ur 
Mathematik (Bonn) and IH\'ES (Bures-sur-Yvette, France)\endthanks
\address Independent University of Moscow, Bolshoi Vlasievskii Pereulok, 11, 
121002 Moscow, Russia\endaddress
\address Steklov Mathematical Institute, Moscow\endaddress
\address Moscow State University, Department of Mathematics\endaddress
\address Present address: Institut des Hautes \'Etudes Scientifiques, 
 35 Route de Chartres, 91440 Bures-sur-Yvette, France\endaddress
\endtopmatter 
\document 
\define\var{\varepsilon}

 \head \S1 Main results and history \endhead
Denote $\Bbb T^n=\Bbb R^n\slash2\pi\Bbb Z^n$.  
Consider a two-dimensional parallel plane foliation on $\Bbb R^n$. 
The standard projection $\Bbb R^n\to\Bbb T^n$ induces 
a foliation on the torus $\Bbb T^n$. This foliation will be denoted by $F$. 

Let $g$ be arbitrary (smooth) Riemann metric on $\Bbb T^n$. It induces 
a complex structure on  each leaf of the foliation $F$ (this 
complex structure  will be  denoted by $\sigma$). The leaves 
are parabolic Riemann surfaces with respect to $\sigma$:  
their universal coverings are conformally equivalent to complex plane.   
This follows from Riemann quasiconformal mapping theorem and boundedness  
of the dilatation of the metric $g$ with respect to the standard  
Euclidean metric of $\Bbb T^n$.  Thus, each leaf admits a complete 
$\sigma$- conformal 
flat metric. This means that there exists a real positive function $\phi$ on 
the leaf such that the metric $\phi g$ on the leaf is flat and complete.  
This function is unique up to multiplication by constant. 

\'E.Ghys \cite{1} stated the following problem: is it true that for any smooth metric $g$ 
on the torus $\Bbb T^3$ the corresponding function $\phi$  
(which defines the flat metric $\phi g$ on the leaves) may be 
chosen to be continuous (smooth) on the whole torus: continuous 
(smooth) not only in the 
parameter of an individual leaf, but also in the transversal parameter 
(in particular, in 
all the intersections of the given leaf with  a  transversal circle)?

In the simplest case, when the leaves of the foliation $F$ are tori, 
the answer to this question is positive. This follows from the classical 
theorem on dependence of uniformization on parameter of complex structure 
\cite{2,3} (Theorem 4 in Subsection 2.4).

We prove the following Theorem giving the positive answer to Ghys' question. 

\proclaim{Theorem 1} Let $F$ be as at the beginning of the paper,  $g$ 
be an analytic /$C^{\infty}$ /measurable uniformly bounded together with 
dilatation along the leaves of $F$/ metric on 
$\Bbb T^n$. There exists an analytic /$C^{\infty}$ /$L_1$/ positive 
function $\phi:\Bbb T^n\to\Bbb R_+$ such that the restriction of the metric 
$\phi g$ to each leaf of the foliation $F$ is flat (in the sense of 
distributions in the third case).
\endproclaim

 Theorem 1 is proved in Section 2.

For a foliation $F$ satisfying appropriate Diophantine  
condition (from one of Definitions 2, 3 or 4) we prove the existence of a 
$C^{\infty}$ (analytic) {\it Euclidean metric on the torus} such that 
the leaves are totally geodesic and its  
restriction to each leaf is conformal and flat (Theorems 2 and 3). This is 
equivalent 
to the statement that the triple consisting of the torus, the foliation $F$ 
and the family of complex structures on its leaves is isomorphic to the analogous triple corresponding 
to a linear foliation on another torus with the complex structure on the leaves 
defined by the standard Euclidean metric. This Diophantine condition is exact.

\definition{Definition 1} We say that a number $\alpha\in\Bbb R\setminus\Bbb Q$ is   
{\it Diophantine}, if there exist constants $C>0$, $s>1$ such that for any pair 
 $m,k\in\Bbb Z$, $k\neq0$, the following inequality holds:
$$|\alpha-\frac mk|>\frac C{|k|^{s+1}}.$$ 
\enddefinition 
 
\definition{Definition 2} Let $x=(x_1,x_2,x_3)$ be coordinates in the 
space $\Bbb R^3$. Consider the foliation on $\Bbb R^3$ 
by level planes of the linear function  
$l(x)=a_1x_1+a_2x_2-x_3$. Let $F$ be the corresponding factorized 
foliation on $\Bbb T^3=\Bbb R^3\slash 2\pi\Bbb Z^3$. We say that $F$ is 
{\it Diophantine}, if there exist constants $C>0$, $s>1$ such that for any 
 $N=(p,m,k)\in\Bbb Z^3\setminus0$ the following inequality holds:
$$|p+ka_1|+|m+ka_2|>\frac C{|N|^s},\ \ |N|=|p|+|m|+|k|.\tag1$$
\enddefinition

\example{Example 1} In the notations of the previous Definition let the 
additive subgroup in $\Bbb R$ generated by 
$a_1$ and $a_2$ contain a Diophantine number.  Then the foliation $F$ is 
Diophantine. 
\endexample

Let us extend the Definition of Diophantine foliation to higher dimensions.

\definition{Definition 3} Let $x=(x_1,\dots,x_n)$ be coordinates in $\Bbb R^n$. 
Consider the foliation on $\Bbb R^n$ by level planes of the linear $n-2$- 
dimensional vector function  
$l(x)=(l_3,\dots,l_n)(x)$, $l_j(x)=a_1^jx_1+a_2^jx_2-x_j$. Let $F$ be 
the corresponding factorized 
foliation on $\Bbb T^n=\Bbb R^n\slash 2\pi\Bbb Z^n$. We say that $F$ is 
{\it Diophantine}, if there exist constants $C>0$, $s>1$ such that for any 
 $m=(m_1,\dots,m_n)\in\Bbb Z^n$ the following inequality holds:
$$|m_1+\sum_{j=3}^na_1^jm_j|+|m_2+\sum_{j=3}^na_2^jm_j|>\frac C{|m|^s},\ \ 
|m|=\sum_{j=1}^n|m_j|.$$
\enddefinition 

\remark{Remark 1} The previous Definition of Diophantine foliation 
coincides with Definition 2 in the case, when $n=3$. 
\endremark

\proclaim{Theorem 2} Let $F$ be a Diophantine foliation (see Definitions 2 and 
3), $g$ be an analytic ($C^{\infty}$) metric on $\Bbb T^n$, $\sigma$ be the  
family of complex structures on the leaves of $F$ such that $g$ is $\sigma$-  
conformal. There exists an analytic ($C^{\infty}$) Euclidean metric on $\Bbb T^n$ such that all the 
leaves are totally geodesic and its restriction to each leaf is $\sigma$- 
conformal. Or equivalently, there exists a discrete rank $n$ additive  
subgroup $G\subset\Bbb R^n$ and an analytic ($C^{\infty}$) diffeomorphism  
$\Bbb T^n\to \Bbb T_G=\Bbb R^n\slash G$ that transforms $F$ to a linear  
foliation and $\sigma$ to the standard complex structure induced by  
the standard Euclidean metric. Inversely, if a linear foliation on $\Bbb T^n$ is 
not 
Diophantine, then there exists a $C^{\infty}$ metric $g$ on the torus such that 
there is no $C^2$ Euclidean metric on $\Bbb T^n$ satisfying the previous 
statement.
\endproclaim 

Theorem 2 is proved in Section 3. 

We give the following weaker Diophantine condition on $F$ necessary and 
sufficient to satisfy the analytic version of the first statement of Theorem 2 
for arbitrary analytic metric.

\definition{Definition 4} In the notations of the previous Definition a 
foliation $F$ is said to be {\it weakly Diophantine}, if 
$$\underline\lim_{|m|\to\infty}(|m_1+\sum_{j=3}^na_1^jm_j|+|m_2+
\sum_{j=3}^na_2^jm_j|)^{\frac1{|m|}}=1.\tag2$$
\enddefinition

\remark{Remark 2} The limit in (2) is always less or equal to 1.
\endremark 

\redefine\tt{\Bbb T^n}
 
 \remark{Remark 3} A Diophantine foliation is always weakly Diophantine.
 \endremark

 \proclaim{Theorem 3} Let $F$ be a weakly Diophantine foliation 
(see Definition 4). 
Then for any analytic metric $g$ on $\Bbb T^n$ there exists an analytic 
Euclidean metric on $\tt$ that satisfies 
the first statement of Theorem 2. Inversely, if $F$ is not 
weakly Diophantine,  then there exists an analytic metric $g$ on $\tt$  
such that there is no $C^2$ Euclidean metric on $\tt$ that satisfies 
 the first statement of Theorem 2. 
 \endproclaim

Theorem 3 is proved in Section 3. It implies the analytic version of Theorem 2.

Earlier a particular case of the problem on existence of continuous complete  
conformal flat  
metric stated at the beginning of the paper was studied by \'E.Ghys \cite{1}.  
He proved the statement of Theorem 1 in the dimension 3 
under the additional assumption that  
in the notations of Definition 2 the additive group generated by the numbers  
$a_1$ and $a_2$ contains either a nonzero rational, or a  
Diophantine number. In the same paper he constructed a counterexample to 
Theorem 1 that is a compact complex manifold equipped with a 
 holomorphic foliation by Riemann surfaces such that for any metric $g$ 
 on the manifold there is no {\it measurable} function $\phi$ such that 
 the restriction to each leaf of the metric $\phi g$ would be flat and 
 complete. Recently he proposed the following general  
 
 \proclaim{Question} Let a compact manifold 
 equipped with a two-dimensional foliation admit at least 
 one smooth metric flat along the leaves. Is it true that for any other 
 smooth metric $g$ 
 there exists a smooth function $\phi$ such that the restriction to the leaves 
 of the metric $\phi g$ is flat?
\endproclaim 

In 1980 A.Haefliger \cite{6} proved a theorem implying the statement analogous 
to Theorem 2: any smooth family of Riemann metrics on the leaves of a linear 
foliation $F$ 
on $\Bbb T^3$ is the restriction to the leaves of a metric on $\Bbb T^3$ for 
which the leaves are minimal, if the 
additive subgroup generated by the numbers $a_j$ (see the notations of 
Definition 2) contains a Diophantine number. 

\head \S2 Existence of fiberwise flat conformal metric. 
Proof of Theorem 1 \endhead
 
 \subhead 2.1. Scheme of the proof of Theorem 1\endsubhead

Theorem 1 is implied by the following more general statement. 

\proclaim{Lemma 1} Let $F$ be at the beginning of the paper, 
$\sigma$ be a family of almost complex structures on the leaves 
of $F$ that is analytic /$C^{\infty}$ /measurable with uniformly bounded 
dilatation/ on $\Bbb T^n$. There exists a nowhere vanishing  
analytic /$C^{\infty}$ /$L_2$/ differential 1- form $\omega$ on 
$\Bbb T^n$ such that its restriction to each leaf is $\sigma$- holomorphic  
and moreover there exists  
a quasiconformal homeomorphism of the universal covering of the leaf 
onto complex plane, whose  
derivative is the pullback of the form $\omega$ under the covering 
(in the sense of distributions in the measurable case).
\endproclaim 

 Lemma 1 is proved in Subsections 2.2 ($C^{\infty}$ case) and 2.3 
(analytic and measurable cases). 

In the conditions of Theorem 1 let $\sigma$ be the family of complex 
structures on the leaves defined by the metric $g$. Let $\omega$ be the 
correspondent $\sigma$- holomorphic 1- form from Lemma 1. 
The metric $\phi g$ from Theorem 1 we are looking for is $\omega
\overline{\omega}$.

The method of the proof of Lemma 1 presented below 
yields a proof of a version of the $C^{\infty}$ Ahlfors-Bers theorem \cite{4} 
(Theorem 4 stated and proved in Subsection 2.4) that is simpler than the 
original one from \cite{4}. 
 \subhead 2.2. Homotopy method. Proof of Lemma 1\endsubhead

For simplicity we present a proof of Lemma 1 only in the case, when $n=3$. 
This proof remains valid in higher dimensions with obvious changes. 

\redefine\tt{\Bbb T^3}

Everywhere below we consider that the leaves of the foliation $F$ under 
consideration are dense in $\tt$, or equivalently, they are not tori. In the 
opposite case the statement of Lemma 1 follows from the classical 
theorem on dependence of uniformization on parameter of complex structure 
\cite{2,3} (cf. Theorem 4 in Subsection 2.4). 

We use the notations of Definition 2 (the coordinates $(x_1,x_2,x_3)$ are chosen 
so that the $x_3$- axis is transversal to the lifting of the foliation 
$F$ to the space). Define the complex coordinate $z=x_1+ix_2$ on the leaves of 
the foliation $F$ and its lifting (the coordinates $x_i$, $z$ and their 
projection pushforwards to the torus will be denoted by the same symbols $x_i$ 
and $z$). 

\define\om{\omega}
\define\omu{\om_{\mu}}
\redefine\c{\Bbb C}
In this Subsection we prove the smooth version of Lemma 1. Its analytic and 
measurable 
versions are proved analogously (the corresponding modifications needed are 
specified in Subsection 2.3). The almost complex structure $\sigma$ on the 
leaves is defined by its $\Bbb C$- linear form that is a complex-valued 1- form 
$\omu=dz+\mu(x)d\bar z$, $\mu\in C^{\infty}(\tt)$, $|\mu|<1$.    
For the proof of Lemma 1 we will construct a $C^{\infty}$ function 
$f:\tt\to\c\setminus0$ such that the restriction to the leaves of the form 
$f\omu$ is closed. Then $\omega=f\omu$ is the $\sigma$- holomorphic differential we are 
looking for. The restriction to each leaf of the metric 
$\omega\overline{\omega}$ is conformal flat 
and complete, as does its pullback to the universal covering of the leaf. 
The quasiconformal homeomorphism from the last statement of 
Lemma 1 is the isometry of the universal covering equipped with the metric 
$\omega\overline{\omega}$ and the complex plane 
equipped with the standard Euclidean metric (this isometry is normalized 
to have the derivative $\omega$). To define this isometry, we fix 
a geodesic $\Gamma$ in the universal covering and consider all the geodesics 
intersecting 
$\Gamma$ orthogonally (at least once). The flatness of the metric 
$\omega\overline{\omega}$ implies that 
no pair of the latters is intersected and no one of them intersects $\Gamma$ more 
than at one point. The isometry we are looking for maps all the geodesics 
under consideration (including $\Gamma$) to straight lines preserving 
lengths of their segments and the orthogonal intersections. This map is 
well-defined and is an isometry (by flatness of  
$\omega\overline{\omega}$) and hence transforms $\sigma$ to the standard 
complex structure. Its derivative is equal to $\omega$ up to 
multiplication by constant with unit module. Indeed the pushforward of the form 
$\omega$ under this isometry is holomorphic (since $\omega$ is $\sigma$- 
holomorphic) and hence, is equal to the standard coordinate differential 
times a holomorphic function of module 1.  This function is constant by 
Liouville's theorem on bounded entire functions. 
Thus, one can achieve that the derivative of the isometry be equal to 
$\omega$ by applying multiplication by appropriate 
complex number with unit module in the image $\Bbb C$.

\define\onu{\omega_{\nu}}

To construct a function $f$ as in the previous item, we use the homotopy method. 
Namely, we include the complex structure $\sigma$ on the leaves 
into the one-parametric 
family of complex structures (denoted by $\sigma_{\nu}$) defined by  their 
$\Bbb C$- linear  
1- forms $\onu=dz+\nu(x,t)d\bar z$, $\nu(x,t)=t\mu(x)$, $t\in[0,1]$. 
The complex structure on the leaves that corresponds to the parameter value $t=0$ is the 
standard one, the given structure $\sigma=\sigma_{\mu}$ corresponds to the 
value $t=1$. 
We will find a $C^{\infty}$ family $f(x,t):\tt\times[0,1]\to\c\setminus0$ of 
complex-valued nowhere vanishing $C^{\infty}$ functions on $\tt$ depending on 
the same 
parameter $t$, $f(x,0)\equiv1$, such that the differential forms $f\onu$ are 
closed on the leaves. Then the function $f(x,1)$ is a one we are looking for. 

We will find a $C^{\infty}$ family $f(x,t)$ of complex-valued 
{\it nonidentically-vanishing} $C^{\infty}$ 
functions on $\tt$ such that the forms $f\onu$ are closed on the leaves. 
Then it follows that the functions $f(x,t)$ vanish 
nowhere. Let us prove this inequality by contradiction. Suppose the contrary. 
Then the set of the parameter values corresponding to the functions $f(x,t)$ 
on the torus having zeroes is nonempty (denote this set by $M$). Its complement 
$[0,1]\setminus M$ is open by definition. Let us show that the set $M$ is open 
as well. This will imply that the parameter segment is a union of two disjoint 
open sets, which will bring us to contradiction. For any 
parameter value the restriction to each leaf of the form $f\onu$ is a 
$\sigma_{\nu}$- holomorphic differential that does not 
vanish identically (by the assumptions that the leaves are dense and for any 
fixed $t$\ $f(x,t)\not\equiv0$ on $\Bbb T^3$). On the other 
hand, each point of the torus has a neighborhood in the corresponding leaf that 
admits a continuous family of $\sigma_{\nu}$- 
holomorphic univalent coordinates (depending on the same 
parameter $t$) by the local $C^{\infty}$ version of the Ahlfors-Bers 
theorem \cite{4} (more precisely, by  
Theorem 4 (or by its weaker version Proposition 2) both stated and proved in 
Subsection 2.4). 
{\it This is the only place in the proof of the $C^{\infty}$ version 
of Lemma 1 we use the Ahlfors-Bers theorem}. The statement that $M$ is open  
follows from the existence of continuous family of 
local $\sigma_{\nu}$- holomorphic 
charts (from Proposition 2), continuity of the family $f\onu$ and 
openness of the set of nonidentically-vanishing holomorphic 
differentials having zeroes. This proves the inequality $f(x,t)\neq0$. 

Let us write down explicitly the fiberwise closeness condition on the form 
$f\onu$. To do this, let us introduce the following

\define\dbz{D_{\bar z}}

\definition{Definition 5} Let $F$ be a linear foliation on $\tt$. In the 
notations of Definition 2 let $z=x_1+ix_2$. 
Let us equip each leaf of the foliation $F$ with the local coordinate $z$.   
Define $D_z$ ($\dbz$) to be the differential operator in the space of 
complex-valued smooth functions on $\tt$ acting as follows: first restrict 
a function to a leaf and then apply the operator 
$\frac{\partial}{\partial z}$ (respectively, $\frac{\partial}{\partial \bar z}$) 
to its restriction to the leaf. 
\enddefinition

\remark{Remark 4} In the conditions of Definition 5
$$D_z=\frac{\partial}{\partial z}+\frac{a_1-ia_2}2\frac{\partial}{\partial x_3},\ \ 
\dbz=\frac{\partial}{\partial\bar z}+\frac{a_1+ia_2}2\frac{\partial}{\partial x_3}\tag3$$
(here $\frac{\partial}{\partial z}$, $\frac{\partial}{\partial\bar z}$ are 
 the operators acting on functions on $\tt$ in the coordinates $(z,x_3)$, 
 not on their restrictions to the leaves). 
\endremark

\remark{Remark 5} In the conditions of Definition 5 let $f,\nu:\tt\to\c$ be 
$C^{\infty}$ functions, $\onu=dz+\nu d\bar z$. The restriction to the leaves 
of the foliation $F$ of the differential form $f\onu$ is closed, iff 
$$\dbz f=D_z(\nu f).\tag4$$
\endremark

Thus, Lemma 1 is implied by the previous discussion and the following 

\proclaim{Lemma 2} Let $F$ be a linear foliation on $\tt$ with dense leaves, 
$\nu(x,t):\tt\times[0,1]\to\c$ be a $C^{\infty}$ family of $C^{\infty}$ functions on $\tt$ 
depending on the parameter $t\in[0,1]$, $\nu(x,0)\equiv0$, $|\nu|<1$. 
There exists a $C^{\infty}$ family 
$f(x,t)$ of complex-valued $C^{\infty}$ functions on $\tt$ depending on the same 
parameter $t$, $f(x,0)\equiv1$, that do not vanish identically and satisfy (4). 
\endproclaim

\demo{Proof} For a family $u(x,t)$ of functions on the torus denote 
$\dot u=\frac{\partial u}{\partial t}$. 
The differentiation of equation (4) in $t$ yields 
$$\dbz\dot f-(D_z\circ\nu)\dot f=(D_z\circ\dot\nu)f,\tag5$$ 
where $D_z\circ\nu$ ($D_z\circ\dot\nu$) is the composition of the operator 
of the multiplication by the function $\nu$ (respectively, $\dot\nu$) and the 
operator $D_z$. 
Any solution $f$ of equation (5) with the initial condition $f(x,0)\equiv1$ 
that does not vanish identically in the torus for no value of $t$ is a one we 
are looking for. 

The proof of Lemma 2 is based on the following properties of the operators $D_z$ 
and $\dbz$. 

\remark{Remark 6} In the conditions of Definition 5 the operators 
$D_z$ and $\dbz$ have the common eigenfunctions $e^{i(N,x)}$, 
$N=(p,m,k)\in\Bbb Z^3$ (which form an orthogonal  
base in the space $L_2(\tt)$). The corresponding eigenvalues are equal to 
$$\lambda_N=\frac i2(p+ka_1-i(m+ka_2)) \ \text{and}\ \lambda'_N=
-\overline{\lambda}_N\tag6$$
respectively. This follows from (3). In particular, $|\lambda_N|=
|\lambda'_N|$ for any $N$. The 
leaves of the foliation $F$ are dense, iff $\lambda_N\neq0$ for all $N\neq0$. 
\endremark

\proclaim{Corollary 1} In the conditions of Definition 5 there exists 
a unitary operator $U: L_2(\tt)\to L_2(\tt)$ preserving averages and 
 the spaces of $C^{\infty}$ (analytic) functions such that 
"$U=\dbz^{-1}\circ D_z$" (more precisely, $U\circ\dbz=\dbz\circ U=D_z$). This 
operator is unique, iff the leaves of the foliation $F$ are dense.
The operator $U$ commutes with partial differentiations and extends up to a 
unitary operator to any Hilbert Sobolev space of functions on $\Bbb T^3$. 
\endproclaim

The operator $U$ from Corollary 1 is defined to have the eigenfunctions 
$e^{i(N,x)}$; the corresponding eigenvalue is equal to 
 $\frac{\lambda_N}{\lambda'_N}= -\frac{\lambda_N}{\overline{\lambda_N}}$  
 (see (6)) if $\lambda_N\neq0$ and to 1 otherwise, i.e., if $D_ze^{i(N,x)}=0$ 
 (in this case it could be taken to be arbitrary number with unit 
 module, if $N\neq0$). 

Let us write down equation (5) in terms of the new operator $U$. Applying 
the "operator" $\dbz^{-1}$ to (5) and substituting $U=\dbz^{-1}\circ D_z$ yields  
$$(Id-U\circ\nu)\dot f=(U\circ\dot\nu)f.$$ 
This equation implies (5). The operator $Id-U\circ\nu$ in the left-hand side 
of this equation is invertible in $L_2(\tt)$ for any $t$ and the norm of the 
inverse operator is bounded uniformly in $t$, since $U$ is unitary and the 
module  
$|\nu|$ is less than 1 and bounded away from 1. As it is shown below (in Proposition 1), so is 
it in all the Hilbert Sobolev spaces $H^j(\Bbb T^3)$. Thus, the last equation 
can be rewritten as 
$$\dot f=(Id-U\circ\nu)^{-1}(U\circ\dot\nu)f,\tag7$$  
 which is an ordinary differential equation in $f\in L_2(\tt)$ with the 
 right-hand side 
 having uniformly bounded derivative in $f$, and so is it in $f$ belonging to 
 arbitrary space $H^j(\Bbb T^3)$ (with respect to the Sobolev scalar product). 
 Therefore, equation (7) written in arbitrary Hilbert Sobolev space has a 
 unique solution 
 with a given initial condition, in particular, with $f(x,0)\equiv1$ (the 
 theorem on existence and uniqueness of solution of ordinary 
 differential equation in Banach space with the right-hand side having 
 uniformly bounded derivative \cite{5}). This solution does not vanish 
 identically in $\tt$  and belongs to 
 all the spaces $H^j(\Bbb T^3)$ for each value of $t$. Therefore it is 
 $C^{\infty}(\tt)$ for any $t$ by Sobolev embedding 
 theorem (see \cite{5}, p.411). Thus, Lemma 2 is implied by the following
 
 \proclaim{Proposition 1} Let $x=(x_1,x_2,x_3)$ be affine coordinates in 
 $\Bbb R^3$, 
 $\tt=\Bbb R^3\slash2\pi\Bbb Z^3$. 
 Let $s\geq0$, $s\in\Bbb Z$, $U$ be a linear operator 
 in the space of $C^{\infty}$ functions on $\tt$ that commutes with the operators 
 $\frac{\partial}{\partial x_i}$, $i=1,2,3$, and extends to any Sobolev space 
 $H^j=H^j(\tt)$, $0\leq j\leq s$, up to a unitary operator. Let $0<\delta<1$, 
 $\nu\in C^s(\tt)$ be a complex-valued function, 
 $|\nu|\leq\delta$. The operator $Id-U\circ\nu$ is invertible and the 
 inverse operator 
 is bounded in all these spaces $H^j$. For any $0<\delta<1$, $j\leq s$, there 
 exists a constant $C>0$ (depending only on $\delta$ and $s$) such that for any 
 complex-valued  function $\nu\in C^s(\tt)$ with $|\nu|\leq\delta$  
 $$||(Id-U\circ\nu)^{-1}||_{H^j}\leq C(1+\max\{|\frac{\partial^k\nu}{\partial x_{i_1},
 \dots,\partial x_{i_k}}|^j\ \ |\ \  \ k\leq j\}).$$
 \endproclaim
 
 \demo{Proof} Let us prove Proposition 1 for $s=1$. For higher $s$ its proof is 
 analogous. 
 
  By definition,  
 $$||U\circ\nu||_{L_2}\leq\delta<1.\tag8$$ Hence, the operator 
 $Id-U\circ\nu$ is invertible in $L_2=H^0$ and 
 $$(Id-U\circ\nu)^{-1}=Id+\sum_{k=1}^{\infty}(U\circ\nu)^k:\tag9$$
 the sum of the $L_2$ operator norms of the sum entries in the right-hand 
 side of (9) is  finite by (8). 
 Let us show that the operator in the right-hand side of (9) is well-defined 
 and bounded in $H^1$. 
 To do this, it suffices to show that the sum of the operator $H^1$- norms of 
 the same entries is finite. 
 
 Let $f\in H^1$. Let us estimate the $H^1$- norm 
 of the images 
 $(U\circ\nu)^kf$. We show that for any $k\in\Bbb N$ 
$$||\frac{\partial}{\partial x_r}((U\circ\nu)^kf)||_{L_2}<ck\delta^{k-1}
||f||_{H^1}, \ \ c=\delta+\max|\frac{\partial\nu}{\partial x_r}|,\ \ 
r=1,2,3.\tag10$$ 
 This will imply the finiteness of the operator $H^1$- norm of the sum in the right-hand 
 side of (9) and Proposition 1 (with $C=4\sum_{k\in\Bbb N}k\delta^{k-1}$).  

Let us prove (10), e.g., for $r=1$. The correspondent derivative in the 
left-hand side of (10) is equal to 
$$(U\circ\nu)^k\frac{\partial f}{\partial x_1}+\sum_{i=1}^k(U\circ\nu)^{k-i}
\circ(U\circ\frac{\partial\nu}{\partial x_1})\circ(U\circ\nu)^{i-1}f\tag11$$
(since $U$ commutes with the partial differentiation by the condition of 
Proposition 1). The $L_2$- norm of the first term in (11) is not greater than 
$\delta^k||f||_{H^1}$ by (8).  That of each entry of the sum in (11) is not greater 
than $\delta^{k-1}\max|\frac{\partial\nu}{\partial x_1}|||f||_{L_2}$ by (8). 
This proves (10). Proposition 1 is 
proved. Lemma 2 is proved, as are the $C^{\infty}$ versions of Lemma 1 and 
Theorem 1.
\enddemo\enddemo

 \define\wt#1{\widetilde#1}
 \define\tdv{\Bbb T^2}

\subhead 2.3. Measurable and analytic cases of Lemma 1 \endsubhead 

In the case, when the metric $g$ is measurable, the proof of Lemma 1 remains 
valid with obvious changes 
(e.g. all the differential equations are understood in the sense of 
distributions) except that of its last statement on existence of quasiconformal 
homeomorphism. Let us prove this statement. 

Let $\omu=dz+\mu d\bar z$ be the 1- form $\c$- linear with respect to $\sigma$, 
$\nu=t\mu$, $t\in[0,1]$. Denote by $f(x,t)$ the solution of ordinary 
 differential equation (7) in $L_2(\tt)$ with $f(x,0)\equiv1$ and put 
$f=f(x)=f(x,1)$. The form 
$f\omu$ is a one we are looking for. Indeed it is closed on the leaves (in the 
sense of distributions on the torus) 
by construction. Let us show that its lifting to the universal covering 
of a generic leaf (with respect to the Lebesgue transversal measure) is the 
derivative of a quasiconformal homeomorphism that 
maps the universal covering onto complex plane. (The universal covering of a leaf 
$L$ will be denoted by $\widetilde L$, the lifting to $\wt L$ of a function $f$ 
(a form $\omu$) on $L$ will be denoted by the same symbol $f$ ($\omu$)). To do 
this, let us approximate the function $\mu$ by functions $\mu_k\in 
C^{\infty}(\tt)$ with moduli less than 1 and uniformly bounded away from 1: 
$|\mu_k|<\delta<1$, $\mu_k\to\mu$ almost everywhere, as $k\to\infty$. Let 
$f_k(x,t)$ be the solution of (7) with $\nu=t\mu_k$, $t\in[0,1]$, 
$f_k(x,0)\equiv1$. Put $f_k(x)=f_k(x,1)$. Then $f_k\to f$  in $L_2(\tt)$: 
the solution of (7) in $L_2(\tt)$ with unit initial condition depends 
continuously on the functional parameter $\nu$ by theorem on dependence of 
solution of ordinary differential equation in Banach space on the parameter 
\cite{5}. Let us fix a 
generic leaf $L$ so that the set $f\neq0$ has a positive measure in $\wt L$ 
and $f_k\to f$ in $L_2$ in compact subsets of $\wt L$. Let us fix a 
point $y_0\in\wt L$. Let $\psi_k:\wt L\to\c$ be the diffeomorphism with 
the derivative $f_k\om_{\mu_k}$ such that $\psi_k(y_0)=0$ (the diffeomorphism 
satisfying these conditions is unique and quasiconformal by the smooth version 
of Lemma 1 proved before). By construction, 
the derivatives $f_k\om_{\mu_k}$ of the diffeomorphisms $\psi_k$ converge to 
$f\omu$ in the sense of distributions on $\wt L$. Let us show that 
the sequence $\psi_k$ converges uniformly in compact sets to a quasiconformal 
homeomorphism. This homeomorphism will be a one we are looking for: 
its derivative will be equal to $f\omu$, as the limit of the derivatives 
$f_k\om_{\mu_k}$.  
 To do this, let us fix another point $y_1\in\wt L$, 
$y_1\neq y_0$. Consider the sequence $b_k=(\psi_k(y_1))^{-1}$ and denote 
$\wt \psi_k=b_k\psi_k$. The maps $\wt \psi_k$ are quasiconformal diffeomorphisms; each of 
them  
transforms the complex structure on $\wt L$ defined by the form $\om_{\mu_k}$ 
to the 
standard one. By definition, they map the points $y_0$ and $y_1$ to 0 and 1 
respectively. They converge uniformly in compact subsets to a quasiconformal 
homeomorphism that transforms the complex structure defined by the form 
$\omu$ to the standard one (by theorem on continuous dependence of the 
uniformizing quasiconformal homeomorphism on the parameter of complex 
structure \cite{2,3,4}). Therefore their derivatives also converge (in the sense of 
distributions) to a nonzero limit. Hence the sequence $b_k$ also converges to 
a nonzero limit, since the contrary would contradict the convergence 
of the derivatives $f_k\om_{\mu_k}$ of the diffeomorphisms $\psi_k$. Therefore the 
initial sequence $\psi_k$ also converges 
uniformly in compact sets to a quasiconformal homeomorphism with the derivative 
$f\omu$. This proves the measurable version of Lemma 1.

The analytic version of Lemma 1 is proved analogously to its $C^{\infty}$ version 
with the following change. We consider equation (7) in the space of analytic 
functions on a fixed closed "annulus" 
$\{x\in\c^3, |\operatorname{Im}x|\leq r\}\slash2\pi\Bbb Z$ 
containing $\tt$. This space is equipped with the scalar 
product 
$$(f,g)=\sum_Nf_N\bar g_Ne^{|N|r}, \ N\in\Bbb Z^3, \ f(x)=\sum_Nf_Ne^{i(N,x)},\ 
g(x)=\sum_Ng_Ne^{i(N,x)}$$ 
(instead of the Sobolev scalar product considered in the smooth case). 
Lemma 1 is proved. The proof of Theorem 1 is completed. 

\subhead 2.4. A proof of the $C^{\infty}$ Ahlfors-Bers theorem\endsubhead

 Let $\Bbb T^2=\Bbb R^2\slash2\pi\Bbb Z^2$, 
$z$ be a complex coordinate in $\Bbb R^2$. 
Recall the following notation: for a complex-valued function $\nu$ on 
$\tdv$ ($\Bbb R^2$) with $|\nu|<1$ by 
$\sigma_{\nu}$ we denote the almost complex structure on $\Bbb T^2$ ($\Bbb R^2$) 
with the $\c$- linear differential $\onu=dz+\nu d\bar z$. 

In this Subsection we give a proof  of the following version of the 
Ahlfors-Bers theorem \cite{4}.

\proclaim{Theorem 4 \cite{2,3}} Let 
$\nu(z,t): \tdv\times[0,1]\to\c$ be a $C^{\infty}$ family of 
$C^{\infty}$ functions on $\tdv$ depending on the parameter $t\in[0,1]$, 
$|\nu|<1$; $\onu=dz+\nu  d\bar z$. There exists a $C^{\infty}$ family 
$f(z,t):\tdv\times[0,1]\to\c\setminus0$ of complex-valued nonvanishing 
$C^{\infty}$ functions on $\tdv$ such that for any fixed 
$t\in[0,1]$ the form $f\onu$ is closed and its standard projection 
pullback to $\Bbb R^2$ is the differential of a  
diffeomorphism of $\Bbb R^2$ onto $\c$.  
\endproclaim

\demo{Proof} Without loss of generality we consider that 
$\nu(z,0)\equiv0$. The proof of Theorem 4 presented below is a modification of 
that of Lemma 1. As in Lemma 2, there exists a $C^{\infty}$ family 
$f(z,t)$ of complex-valued $C^{\infty}$ functions on $\Bbb T^2$ such that 
$f(z,0)\equiv1$, the forms $f\onu$ are closed 
and for any $t\in[0,1]$\ $f(z,t)\not\equiv0$. It is the unique solution of the 
ordinary differential equation (7) in $L_2(\Bbb T^3)$ 
 with the initial condition $f(z,0)\equiv1$, where  the operator 
 $U=(\frac{\partial}{\partial\bar z})^{-1}\circ\frac{\partial}{\partial z}$ 
 is defined to have the eigenbase 
 $e^{i(n_1\operatorname{Re} z +n_2\operatorname{Im} z)}$, $n_1,n_2\in\Bbb Z$,  
 with the eigenvalues equal to $\frac{n_1-in_2}{n_1+in_2}$ for 
 $(n_1,n_2)\neq0$ and 
 the unit eigenvalue at the constant function. (This base is formed by the 
 common eigenfunctions of the operators 
 $\frac{\partial}{\partial\bar z}$ and $\frac{\partial}{\partial z}$.) The 
 functions $f(z,t)$ are $C^{\infty}$,   
 since for any $s\in\Bbb N$ and $t\in[0,1]$ $f(z,t)\in H^s(\Bbb T^2)$ and 
 by Sobolev embedding theorem (\cite{5}, p.411). The inclusions $f\in H^s$ 
 follow from boundedness in any Hilbert Sobolev norm of the operator in the 
 right-hand side of (7) (cf. the discussion 
 following (7) and the inequality of Proposition 1, which remains valid 
 for the new operator $U$). Equation (7) 
 implies that the forms $f\onu$ are closed, as in the proof of Lemma 2. 
 
 Let us prove that 
 the family $f$ is a one we are looking for. As at the beginning of Subsection 
 2.2, it suffices to show that $f(z,t)$ vanishes nowhere.  
 The proof of this inequality repeats the one from the proof of Lemma 1 
 and is done by contradiction. Suppose 
 the contrary, i.e., the set of parameter values corresponding to the 
 functions having zeroes (denote it by $M$) is nonempty. Its complement 
 is open. It suffices to show 
 that $M$ is open: this will imply that the parameter segment is a union of two 
 disjoint nonempty open sets, which brings us to contradiction. Zeroes of 
 a function $f(z,t)$ 
 are exactly zeroes of the 1- form $f\onu$, which is $\sigma_{\nu}$- 
 holomorphic and does not vanish identically. 
 The openness of the set $M$ is implied by the continuity of the family $f\onu$ 
 of 1- forms, the openness of the set of 
nonidentically-vanishing holomorphic differentials having zeroes and the 
following local version of the Ahlfors-Bers theorem \cite{4}:
 
 \proclaim{Proposition 2 \cite{4}} Let $\nu(z,t)$ be a $C^{\infty}$ family of 
 complex-valued $C^{\infty}$ functions 
 in unit disc (equipped with the coordinate $z$) depending on the parameter 
 $t\in[0,1]$, $|\nu|<1$. For any $t_0\in[0,1]$ 
 there exist a smaller disc $V$ centered at 0  
 and  a $C^{\infty}$ family of $C^{\infty}$ 
 $\sigma_{\nu}$- holomorphic univalent coordinates in $V$ defined for all $t$ 
 close enough to $t_0$.   
\endproclaim

 Proposition 2 will be proved below, more precisely, it will be 
 reduced to the following statement.

\proclaim{Proposition 3} 
For any $\var>0$ there exists a $\delta>0$ such that 
for any complex-valued  function $\mu\in C^3(\tdv)$ with $||\mu||_{C^3}<\delta$ 
there exists a complex-valued function $f\in C^1(\tdv)$ with 
$||f-1||_{C^1}<\var$ such that the form $f\omu=f(dz+\mu d\bar z)$ is closed 
(i.e., if $\var<1$, then locally it is the differential of a 
$\sigma_{\mu}$- holomorphic univalent function). The  function $f$ 
is unique up to multiplication by constant and can be normalized to 
depend continuously on $\mu\in C^3(\tdv)$ as a $C^1$- valued functional. If 
$\mu$ is $C^{\infty}$, then so is $f$. 
\endproclaim   
 
 \demo{Proof} Firstly let us prove  Proposition 3 with the change of the 
 $C^1$- norm in its statement to the $H^3$- norm. To do this, let us include 
 $\mu$ into the family 
 $\nu(z,t)=t\mu(z)$ of functions depending on the parameter $t\in[0,1]$. 
 As in the proof of Lemma 2, there exists a unique solution $f(z,t)$ of the 
 ordinary differential equation (7) in $H^3(\tdv)$ 
 (with the new operator $U$ defined at the beginning of 
 the Subsection) such that $f(z,0)\equiv1$. The function 
 $f(z)=f(z,1)$ is the one we are looking for. Indeed, the form $f\omu$ 
 is closed,  $f$ is $C^{\infty}$ if so is $\mu$ (as in the proof of 
 Lemma 2). The function $f(z)$ depends 
 continuously on the functional parameter $\mu\in C^3(\tdv)$  as an 
 $H^3$- valued functional (by the inequality of Proposition 1 (which remains 
 valid for the new operator $U$) and the theorem 
 on continuous dependence of solution of bounded ordinary 
 differential equation in Banach space  on the parameter \cite{5}). Therefore 
 $f(z)$ is $H^3$- close to 1, whenever 
 the norm $||\mu||_{C^3}$ is small enough. This together with continuity 
 of the canonical embedding $H^3\to C^1$ (Sobolev embedding theorem, \cite{5}, 
 p.411) proves Proposition 3.\enddemo     
 
Let us prove Proposition 2. Without loss of generality we 
consider that $\nu(0,t_0)=0$ (one can achieve this by applying appropriate 
$\Bbb R$- linear 
change of the coordinate $z$). Let us choose arbitrary $\delta'>0$. Without loss 
of generality we assume that all the derivatives of the function 
$\nu(z,t_0)$ in $z$ up to the order 3 have moduli less than $\delta'$ in the disc 
$|z|\leq\frac12$  (one can achieve this by applying appropriate homothety in $z$). 
Let $\delta,\var<1$ be as in Proposition 3. In the case, when $\delta'$ is 
chosen to be small enough, the restriction to the disc $|z|\leq\frac14$ of the 
function $\nu$ 
extends to $\c=\Bbb R^2$ up to a double $2\pi$- periodic function whose 
pushforward to $\tdv$ has $C^3$- norm less than $\delta$. This together with 
Proposition 3 implies 
existence of local family of $\sigma_{\nu}$- holomorphic coordinates we are 
looking for. This proves Proposition 2 and Theorem 4. \enddemo

\head \S3. Diophantine foliations. Proofs of Theorems 2 and 3\endhead
We present only the proof of the three-dimensional $C^{\infty}$ version of 
Theorem 2. This proof remains valid in higher dimensions with obvious changes. 
The proof of 
Theorem 3 (which implies the analytic version of Theorem 2) is analogous with 
the modifications similar to those needed in the proof of the analytic version 
of Lemma 1 (see the end of Subsection 2.3). 

Let $F$ be a Diophantine foliation, $l(x)=a_1x_1+a_2x_2-x_3$ be the 
corresponding linear function 
(see the notations of Definition 2), $z=x_1+x_2$. Let $\sigma$ be a 
$C^{\infty}(\Bbb T^3)$ family of almost complex structures on the leaves of 
$F$ (e.g., defined by a smooth metric on $\Bbb T^3$). Let $\omega$ be the correspondent differential form  
from Lemma 1. It is uniquely defined modulo $dl$ 
up to multiplication by constant, since the leaves 
of the foliation $F$ are  dense ($F$ is Diophantine). Everyone of the two 
first equivalent statements in Theorem 2 is equivalent to the possibility to 
choose the form $\omega$ to be closed not only on the leaves, but on the whole 
torus. Indeed the second one of these statements in Theorem 2 implies the existence of a 
differential 1- form closed on $\Bbb T^3$ and holomorphic on the leaves. 
Inversely, for a given closed fiberwise holomorphic 1- form $\omega$  
the Euclidean metric $\omega\overline{\omega}+dld\bar l$ on 
$\tt$ satisfies the first statement of Theorem 2. Let $\omega=f(x)(dz+\mu 
d\bar z)$ be a fixed $C^{\infty}$ differential form on $\tt$ $\sigma$- holomorphic on 
the leaves. Let us show that there exists a $C^{\infty}$ function 
$h:\tt\to\c$ such that the form 
$$f(dz+\mu d\bar z)-hdl\tag12$$ 
is closed on $\tt$. This will prove the $C^{\infty}$ version of Theorem 2. 
To do this, we use the following equivalent reformulation of the condition of 
closeness of the form (12).

\remark{Remark 7} In the conditions of Definition 5 let 
$f$, $\mu$, $h$ be functions on $\tt$ such that the form (12) is 
closed on the leaves of $F$. Then (12) is closed on $\tt$, iff $f$, $\mu$ 
and $h$ satisfy the following system of differential equations: 
$$\cases
& \frac{\partial f}{\partial x_3}=D_zh\\
& \frac{\partial(\mu f)}{\partial x_3}=\dbz h\endcases.\tag13$$ 
Indeed, let us write down the closeness condition on the form (12) in the 
coordinates $(x_1,x_2,l)$. In the new coordinates $\frac{\partial}{\partial z}=
D_z$, $\frac{\partial}{\partial\bar z}=\dbz$, 
the operator $\frac{\partial}{\partial l}$ coincides with the operator 
$-\frac{\partial}{\partial x_3}$ correspondent to the coordinates $(x_1,x_2,x_3)$. 
By definition, (12) is closed on the leaves (or equivalently, $f$ and $\nu=\mu$ 
satisfy (4)). Under this assumption, the condition that (12) is closed on 
$\tt$ is equivalent to the system of differential equations 
$\frac{\partial f}{\partial l}=-\frac{\partial h}{\partial z}$, 
$\frac{\partial(\mu f)}{\partial l}=-\frac{\partial h}{\partial\bar z}$. 
Rewriting these equations in the coordinates $(x_1,x_2,x_3)$ yields (13). 
\endremark

Let us show the existence of a $C^{\infty}$ solution $h$ to (13). This together 
with the previous Remark will 
prove the first statement of the $C^{\infty}$ version of Theorem 2. To do this, 
we will use the following characterization of $C^{\infty}$ functions on the torus.

\remark{Remark 8} A function $h\in L_1(\tt)$ with the Fourier series 
$\sum_Nh_Ne^{i(N,x)}$ is $C^{\infty}$, iff 
$$\sum_N|N|^s|h_N|<\infty \ \ \text{for any}\ \ s\in\Bbb N.\tag14$$
\endremark

We use the following equivalent reformulation of the Diophantiness condition. 

\remark{Remark 9} Let $F$ be a linear foliation on $\tt$, $D_z$ be the corresponding 
differential operator from Definition 5, $\lambda_N$ be its eigenvalues 
from (6). The foliation $F$ is Diophantine, iff there exist $c>0$, $s>1$  
such that for any $N\in\Bbb Z^3\setminus0$  
$$|\lambda_N^{-1}|\leq c|N|^s.\tag15$$
\endremark

Define the "inverse" operator $D_z^{-1}$ in the space of infinitely-smooth 
functions on $\tt$ with zero 
average (this space is contained in the Hilbert subpace in $L_2$ 
generated by the eigenfunctions $e^{i(N,x)}$ 
of the operator $D_z$ with $N\neq0$) to have the same eigenfunctions $e^{i(N,x)}$ 
with the eigenvalues inverse to those of $D_z$. 
This operator is well-defined at least on the 
eigenfunctions, since the corresponding eigenvalues $\lambda_N$ of the 
operator $D_z$ do not vanish (by (15)). It is well-defined on the space of 
infinitely-smooth 
functions on $\tt$ with zero average by Remark 8 and (15) and has zero kernel. 
The function $h=D_z^{-1}(\frac{\partial f}{\partial x_3})$ is the one we are 
looking for. Indeed it is 
a $C^{\infty}$ function satisfying the first equation in  (13). It satisfies the 
second equation as well. Indeed, applying the operator $D_z$ to 
the second equation and substituting 
$h=D_z^{-1}(\frac{\partial f}{\partial x_3})$ yields (4) with $\nu=\mu$, which 
is satisfied by the assumption that the form $f\omu$ is closed on the leaves. 
Therefore, the $D_z$ images of both parts of the 
second equation coincide. Hence, these parts coincide theirselves: 
 they have zero average,  and hence they are obtained from their (coinciding) 
 $D_z$ images by applying the operator $D_z^{-1}$. 
 This proves the first statement of the $C^{\infty}$ version of Theorem 2.

Now let us prove the last statement of Theorem 2. Let 
$F$ be a linear nondiophantine foliation on $\tt$. Let us prove the existence 
of a $C^{\infty}(\Bbb T^3)$ family of almost complex structures on the leaves of 
the foliation $F$  
such that there is no $C^2$ Euclidean metric on $\tt$ satisfying the first statement 
of Theorem 2 (this is equivalent to the last statement of Theorem 2). 
In the case, when the leaves are not dense 
(i.e., they are tori), for a generic $C^{\infty}(\Bbb T^3)$ family of complex 
structures  
the leaves do not have the same conformal type. This is the family we are 
looking for. Indeed there is no diffeomorphism satisfying the second one 
of the equivalent statements in Theorem 2, since otherwise all the leaves would 
be conformally equivalent to each other. 

Everywhere below we consider that the foliation $F$ has dense leaves. We 
consider that its lifting 
to the 3-space is transversal to the $x_3$- axis and we use the notations of 
Definition 2. In the proof of the last statement of Theorem 2 we use the following 

\remark{Remark 10} Let $F$ be a linear foliation on $\tt$ with dense leaves. 
In the notations of Definition 2 let $z=x_1+ix_2$, $\mu\in C^{\infty}(\tt)$, 
$|\mu|<1$, $\omu=dz+\mu d\bar z$. There exists a unique complex-valued 
function $f\in C^1(\tt)$ 
such that the restriction to the leaves of the form $f\omu$ 
is closed, up to multiplication by constant. Indeed, there exists at least one 
function $f\in C^{\infty}(\tt)$ satisfying this 
condition and nowhere vanishing (this is an equivalent reformulation of the first 
statement of Lemma 1). Let $\varphi\in C^1(\tt)$ be another function such that the form 
$\varphi\omu$ is closed on the leaves. 
 Then the ratio $\frac{\varphi}f$ is $\sigma_{\mu}$- holomorphic and  
uniformly bounded on each leaf. Therefore, it is constant on each leaf (the leaves are 
parabolic Riemann surfaces), and hence, on $\tt$ as well (by density of the 
leaves). \endremark

For the proof of the last statement of Theorem 2 in the $C^{\infty}$ case we show the 
existence of complex-valued functions $\nu,f\in C^{\infty}(\tt)$, 
$|\nu|<1$, such that the form $f\onu$ is closed on the leaves (or equivalently, 
$f$ and $\nu$ satisfy (4)) and that there is no $C^2$ (and even $L_2$) 
complex-valued 
function $h$ on $\tt$ such that the correspondent form (12) with $\mu=\nu$  
is closed on $\Bbb T^3$ (or equivalently, (13) holds, cf. Remark 7).  
This will prove the $C^{\infty}$ version of Theorem 2. 

Let $\lambda_N$ be the eigenvalues of the operator $D_z$ (see (6)). To 
construct a pair $(\nu,f)$ as in the previous paragraph, let us choose 
and fix a sequence of distinct multiindices $N_j=(p_j,m_j,k_j)\in\Bbb Z^3$ with 
$k_j\neq0$ such that for any $s\in\Bbb N$ 
there exists a number $J\in\Bbb N$ such that for any $j>J$ 
$$|\lambda_{N_j}^{-1}|>|N_j|^s.\tag16$$
The existence of a sequence $N_j$ satisfying the previous estimate (16) 
follows from nondiophantiness of the 
foliation $F$ (which is contrary to (15)). This sequence can be chosen so 
that $k_j\neq0$ for all $j$: one can achieve this by choosing a subsequence 
of multiindices $N_j\neq0$ such that  $|\lambda_{N_j}^{-1}|>1$. 
Then $k_j\neq0$, which is implied by the following inequality: let $N=(p,m,k)\in\Bbb 
Z^3\setminus0$ be such that $|\lambda_N|<1$; then $k\neq0$. This 
inequality follows from (6).   

Consider the family 
$f(x,t)=1+t\sum_{j=1}^{\infty}\frac{\lambda_{N_j}}{k_j}e^{i(N_j,x)}$ 
of functions on $\tt$ depending on real parameter $t$.  This is a $C^{\infty}$ 
family of $C^{\infty}$ functions, since the coefficients of the  Fourier series 
in its formula satisfy (14) (by (16)). As it is shown below, any function 
$f=f(x,t)$ 
corresponding to small enough nonzero value of $t$ is a one we are looking for. 

The statement on the existence of a function $\nu$ satisfying  (4) is proved 
below (Lemma 3). For any fixed nonzero 
value of the parameter $t$ and any complex-valued function $\mu=\nu$ there is 
no function 
$h\in L_2(\Bbb T^3)$ such that the correspondent form (12) is closed, since  
there is no $h\in L_2(\Bbb T^3)$ satisfying the first equation in (13) and 
by Remark 7. Indeed, if an $L_2$ function $h$ 
satisfies the first equation in (13), then its Fourier coefficients with the 
indices $N_j$  would be 
equal to $it$, which is impossible if $t\neq0$. Thus, now Theorem 2 is implied 
by the following 

\proclaim{Lemma 3} Let $F$ be a linear foliation on $\tt$, $f(x,t)$ be a 
$C^{\infty}$ 
family of complex-valued $C^{\infty}$ functions on $\tt$ uniformly bounded away 
from zero depending on real parameter $t$, 
$f(x,0)\equiv1$. There exists a $C^{\infty}$ family $\nu(x,t)$ of complex-valued 
$C^{\infty}$ functions on $\tt$ 
depending on the same parameter $t$, $\nu(x,0)\equiv0$, and satisfying (4).  
\endproclaim

\demo{Proof} Lemma 3 is proved analogously to Lemma 2. Differentiating  
equation (4) in $t$ yields (5). The solution $\nu(x,t)$ of (5) with 
zero initial condition will be a solution of (4). Let $U$ be the operator 
from Corollary 1. Applying subsequently the "operator" $D_z^{-1}$ and the 
multiplication by 
$f^{-1}$ to (5) and substituting $U^{-1}=D_z^{-1}\circ\dbz$ yields 
$$\dot\nu=\frac{U^{-1}\dot f-\nu\dot f}f.$$
The last equation implies (5) and has a unique infinitely-smooth solution 
$\nu(x,t)$ with any  
given $C^{\infty}$ initial condition (e.g., $\nu(x,0)\equiv0$): its right-hand 
side has bounded derivative in $\nu$ in any Hilbert Sobolev norm (cf. the 
proof of Lemma 2).  This proves 
Lemma 3. The proof of the $C^{\infty}$ version of Theorem 2 is completed.
\enddemo


\head Acknowledgements\endhead

I am grateful to Yu.S.Ilyashenko and \'E.Ghys, who have attracted my 
attention to the problem. I am grateful to them and to  E.Giroux, J.Milnor for 
helpful discussions. 
The paper was written while I was staying at the Independent 
University of Moscow and Institut des Hautes \'Etudes Scientifiques 
(IH\'ES, Bures-sur-Yvette, France). My stay at these Institutions was 
supported by the European Post-Doctoral Institute fellowship of 
IH\'ES. I wish to thank both Institutions for the hospitality and 
support.

\head References\endhead 

1. Ghys, \'E. Sur l'uniformisation des laminations paraboliques. - in  
Integrable systems and foliations, ed. C.Albert, R.Brouzet, J.-P. Dufour  
(Montpellier, 1995), Progress in Math. 145 (1996), 73-91. 
 
2. Abikhow, W. Real analytic theory of Teichm\"uller space, - Lect. Notes in  
Math., 820, Springer-Verlag (1980). 
 
3. Ahlfors, L. Lectures on quasiconformal mappings, - Wadsworth (1987).  

4. Ahlfors, L.; Bers, L. Riemann's mapping theorem for variable metrics, - 
Ann. of Math. (2) 72 (1960), 385-404. 

5. Choquet-Bruhat, Y., \ de Witt-Morette, C., \ Dillard-Bleick, M. Analysis, 
Manifolds and Physics, - North-Holland, 1977.

6. Haefliger, A. Some remarks on foliations with minimal leaves. - J. Differential 
Geometry, 15 (1980), no 2, 269-284 (1981). 
\enddocument